\documentclass{birkjour}

\usepackage{enumitem}
\newtheorem{axm}{Axiom}
\newtheorem{thm}{Theorem}[section]

\newtheorem{lem}[thm]{Lemma}

\newtheorem{cnj}[thm]{Conjecture}
\theoremstyle{definition}

\newtheorem{defn}[thm]{Definition}
\theoremstyle{remark}

\numberwithin{equation}{section}
\newcommand{\AUTHOR}{}

\begin{document}

%-------------------------------------------------------------------------
% editorial commands: to be inserted by the editorial office
%
%\firstpage{1} \volume{228} \Copyrightyear{2004} \DOI{003-0001}
%
%
%\seriesextra{Just an add-on}
%\seriesextraline{This is the Concrete Title of this Book\br H.E. R and S.T.C. W, Eds.}
%
% for journals:
%
%\firstpage{1}
%\issuenumber{1}
%\Volumeandyear{1 (2004)}
%\Copyrightyear{2004}
%\DOI{003-xxxx-y}
%\Signet
%\commby{inhouse}
%\submitted{March 14, 2003}
%\received{March 16, 2000}
%\revised{June 1, 2000}
%\accepted{July 22, 2000}
%
%
%
%---------------------------------------------------------------------------
%Insert here the title, affiliations and abstract:
%

\title[Is The Parallel Postulate Necessary For Euclidean Geometry?]
 {Is The Parallel Postulate Necessary For Euclidean Geometry?}

\ifdefined\AUTHOR
%----------Author 1
\author[Chengpu Wang]{Chengpu Wang*}
\email{Chengpu@gmail.com}

%----------Author 2
\author{Alice Wang}
\email{atwang@uchicago.edu}
\fi

\thanks{}

%----------classification, keywords, date

%\keywords{Euclidean, formal system}

\date{Jan 28, 2023}
%----------additions
%\dedicatory{To the students who are willing to reinvent the wheel}
%%% ----------------------------------------------------------------------

\maketitle

\begin{abstract}
As a much later addition to the original Euclidean geometry, the parallel postulate distinguishes non-Euclidean geometries from Euclidean geometry.  
This paper will show that the parallel postulate is unnecessary because the 4th Euclidean axiom can already achieve the same goal. 
%In fact, the 4th Euclidean axiom itself can also be reduced to a theorem.
Furthermore, using the 4th Euclidean axiom can measure space curvature locally on manifold, while using the parallel postulate cannot.
\end{abstract}

\section{Introduction}

Euclidean geometry \cite{Euclid}\cite{Plane}\cite{Old and New} is the origin for modern mathematics. The original Euclidean axiom set contains 4 axioms. 
It is regarded as insufficient \cite{Old and New}\cite{Hilbert}, so it has been enhanced as the Hilbert axiom set \cite{Hilbert} of 21 axioms in 5 groups. 
For common usage, a parallel postulate is added to the Euclidean axiom set as the 5th axiom \cite{Euclid}\cite{Plane}\cite{Old and New}. 
In this paper, the Euclidean axiom set does not contain the 5th axiom.

One conceptual difficulty of Euclidean geometry is the measurement of angle: an angle is measured by the corresponding arc length, while the arc length is calculated by the angle measurement \cite{Plane}, which results in circular definitions \cite{Formal System}. 
Perhaps to overcome this difficulty, later the Birkhoff axiom set \cite{Birkhoff} requires both distance and angle to be independently measurable, which is different from the common practice. 
Also, the Birkhoff axiom set is too simple to be mathematically rigorous.

The Hilbert axiom set is of second-order logic.  
It has an equivalent first-order logic axiom set as the Tarski axiom set \cite{Tarski} which is quite different from the common understanding of Euclidean geometry, so it is not commonly taught in geometry classes. 

All modern axiom sets use the parallel postulate to distinguish non-Euclidean geometries \cite{Non-Euclidean Geometry} from Euclidean geometry.  
However, since the parallel postulate is based on straight lines which are infinitive in extend, it cannot measure the surface curvature \cite{Differential Geometry} locally on a manifold \cite{Differential Geometry}.  
This paper will show that the 4th Euclidean axiom “All right angles are equal to each other” can distinguish non-Euclidean geometries from Euclidean geometry both locally and globally.  Furthermore, it leads to a new and simpler surface curvature analytic called right ratio, to improve over the explicit Gaussian curvature, the implicit Riemann curvature tensor \cite{Differential Geometry}, or other existing analytics.   

Both the Hilbert and the Tarski axiom sets use congruence and similarity in the axiom set, while the original Euclidean axiom set does not but instead concludes congruence and similarity as theorems. 
This paper presents a modernized Euclidean axiom set.

\section{The Modernized Euclidean Axiom Set}

\subsection{Notation}

In this new axiom set, point is the only primitive notion [6]. 
Let a capitalized letter denote a point. Define a geometric object as a set of points. 
A point can be in $\in$ or not in $\not \in$ a geometric object. 
Two geometric objects can either intersect $\cap$ or union $\cup$ each other. 
$\oslash$ is the empty set. 
let $\forall$ for any or all, $\exists$ for some, and $\exists !$ for only one.  

Two mathematical statements can be combined logically to form a new statement: $\bigwedge$ for and, $\bigvee$ for or. 
Two mathematical statements can also have cause-effect relationships: $\Longrightarrow$ for leading to, $\iff$ for equivalent (which is leading to in both ways), and : for common condition or common definition.

\subsection{Measurable}

As the embodiment of the 3rd Euclidean axiom “To describe a circle with any center and distance” \cite{Euclid}:

\begin{defn}[distance]
Between two points $A$ and $B$, the distance $|AB|$ is a non-negative real value which satisfies the following requirements:
\begin{itemize}
\item $\forall A: |AA| = 0$.
\item $\forall A, B: 0 \leq |AB| \bigwedge |AB| = |BA|$.
\item $\forall A, B, C: |AB| \leq |AC| + |CB|$.
\end{itemize}
\end{defn}

\begin{axm}[measurable]
$\forall A, B \Longrightarrow \exists ! |AB|$.
\end{axm}

Axiom (measurable) is equivalent to “Postulate I. Space is metric” \cite{Isometric}, in which a measurable space is a metric space. It is implicit in the Hilbert axiom set \cite{Hilbert} and the Tarski axiom set \cite{Tarski} by asserting “congruence” relation between any two points \cite{Old and New}. 

\begin{defn}[round]
A \emph{round} geometric object is $\underline{A:r} \equiv  \{P: |AP| = r \} $ in which the point $A$ is the \emph{center}, while the constant $r$ is the \emph{radius}.
\end{defn}

\subsection{Continuous}

As the embodiment of the 1st Euclidean axiom “To draw a straight line from any point to any point.” \cite{Euclid}:

\begin{defn}[straight segment]
A straight segment $\underline{AB}$ between two points $A$ and $B$ is defined as:

$\underline{AB} \equiv \{P: |AB| = |AP| + |PB|, |AP| \in [0, |AB|], P \iff |AP| \}$.
\end{defn}

\begin{axm}[continuous]
\label{continuous}
$\forall A, B \Longrightarrow \exists ! \underline{AB}$.
\end{axm}

\subsection{Boundless}

As the embodiment of the 2nd Euclidean axiom “To extend a finite straight line continuously in a straight line.” \cite{Euclid}:
\begin{axm}[boundless]
\label{boundless}
$\forall \underline{AB} \Longrightarrow \exists C, D: A \in \underline{CB} \bigwedge B \in \underline{AD}$.
\end{axm}

\begin{defn}[ray]
A ray $\underline{AB*}$  from a point $A$ passing a different point $B$ is defined as 
$\underline{AB*} \equiv \underline{AB} \cup \{P: |AP| = |AB| + |BP|, P \iff |AP| \}$.
\end{defn}

\begin{defn}[straight line]
A straight line $\underline{*AB*}$ passing two different point $A$ and $B$ is defined as:
$\underline{*AB*} \equiv \underline{*AB} \cup \underline{AB*}$.
\end{defn}

\subsection{Topologically Simplest}

A Euclidean surface is simpler when compared with topologically \cite{Topology} more complex surfaces such as Mobius ring and Klein bottle \cite{Non-Euclidean Geometry}:

\begin{defn}[opposite side]
If $C, D \not \in \underline{*AB*} \bigwedge \underline{*CD*} \cap \underline{*AB*} \neq \oslash$, $C$ and $D$ are on the \emph{opposite side}s of $\underline{*AB*}$, which is denoted as $C, D \div \underline{*AB*}$.
\end{defn}

\begin{defn}[2-dimensional space]
When $C, D \div \underline{*AB*}$, the \emph{2-dimensional space} $\{ C, D \div \underline{*AB*} \} \equiv \underline{*AB*} \cup \{P: P, C \div \underline{*AB*} \} \cup \{P: P, D \div \underline{*AB*} \}$
\end{defn}

\begin{defn}[2-sided 2-dimensional space: surface]
If $\{ C, D \div \underline{*AB*} \}$ further satisfies $\{P: P, C \div \underline{*AB*} \} \cap \{P: P, D \div \underline{*AB*} \} = \oslash$, it is a \emph{2-sided} 2-dimensional space $\{ C, D \div \underline{*AB*} \}_2$, which is simplified as a \emph{surface}.
\end{defn}

\begin{thm}[same side]
$\{ C, D \div \underline{*AB*} \}_2: C, E \div \underline{*AB*} \iff \underline{*DE*} \cap \underline{*AB*} = \oslash$, which is defined as $D$ and $E$ being on the \emph{same side} of $\underline{*AB*}$.
\end{thm}

\subsection{Strictly Monotonic}

\begin{figure}
\includegraphics{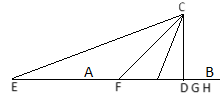}
\caption{Demonstration of Axiom (strictly monotonic).}
\label{fig:strictly_monotonic}
\end{figure}

A straight line intersects a circle in two points on a Euclidean surface.
A straight line is a geodesic \cite{Non-Euclidean Geometry}\cite{Differential Geometry} on a Euclidean surface.
On a non-Euclidean surface, a geodesic may intersect a circle for more than two points.
To capture such simplicity of a Euclidean surface, as shown in Figure \ref{fig:strictly_monotonic}:

\begin{axm}[strictly monotonic]
\label{tangent}
$\forall C \not \in \underline{*AB*}$:
\begin{enumerate}
\item $\exists ! D \Longrightarrow \{D\} = \underline{C:|CD|} \cap \underline{*AB*}$ in which $D$ is the \emph{tangent} between $C$ and $\underline{*AB*}$, which is denoted as $D = C \perp \underline{*AB*}$;

\item $\forall E \in \underline{*AB*}: E \neq D \iff |DE| < |CE|$.

\item $\forall E, F \in \underline{*AB*}: |DF| < |DE| \iff |CF| < |CE|$.
\end{enumerate}
\end{axm}

Axiom (strictly monotonic) establishes perpendicular relationships \cite{Euclid}\cite{Plane} between the two straight lines as $D = \underline{*CD*} \perp \underline{*AB*} = \underline{*AB*} \perp \underline{*CD*}$:

\begin{thm}[perpendicular]
$D = C \perp \underline{*AB*}$:
\begin{itemize}
\item $\forall E \in \underline{*CD*} \Longrightarrow D = E \perp \underline{*AB*}$;

\item $\forall F \in \underline{*AB*} \Longrightarrow D = F \perp \underline{*CD*}$;
\end{itemize}
\end{thm}

\subsection{Angle Measurement}

\begin{defn}[angle]
Any two rays from one point form an angle, in which the point is the \emph{vertex}, while the two rays are the two \emph{arm}s of the angle. 
Letting the vertex be $A$ and the two arms be $\underline{AB*}$ and $\underline{AC*}$, the angle is denoted as $\angle BAC$ or equivalently $\angle CAB$.
\end{defn}

\begin{defn}[angle distance ratio]
To measure an angle $\angle BAC$:
\begin{enumerate}
\item Use a positive value called angle measuring distance $r$ to construct $\underline{A:r}$, which intersects $\underline{AB*}$ and $\underline{AC*}$ at $D$ and $E$ respectively; 

\item The \emph{angle distance ratio} is: $|\angle BAC|_r \equiv |DE|^2/(2r)^2 \in [0,1]$, with $0$ for a ray and $1$ for a straight line.
\end{enumerate}
\end{defn}

\begin{defn}[right angle]
The four angles formed by two perpendicular straight lines are all \emph{right angle}s, such as $\angle ADC$ and $\angle BDC$ for $D = C \perp \underline{*AB*}$.
\end{defn}

\begin{defn}[right ratio]
The angle distance ratio of a right angle is a \emph{right ratio}.
\end{defn}

\begin{defn}[local right ratio]
The limit of the right ratio when the measuring distance $\rightarrow 0$ is the \emph{local right ratio}.
\end{defn}

\subsection{Flat}

As the embodiment of the 4th Euclidean axiom “All right angles are equal to one another” \cite{Euclid}:
\begin{thm}[flat]
If on a surface, all right ratios are a constant, they are $1/2$.
\end{thm}

\begin{defn}[flat]
A \emph{Euclidean surface} is a surface on which all the right ratios are $1/2$.
\end{defn}

\begin{figure}
\includegraphics{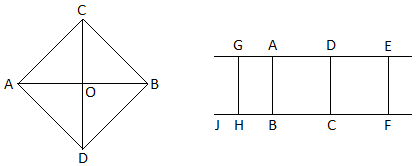}
\caption{ \textsl{Left:} Form a square in a Euclidean surface. \textsl{Right:} Stack identical squares side by side to form rectangles }
\label{fig:squares}
\end{figure}

In Figure \ref{fig:squares} \textsl{Left}, let $O = \overline{AB} \perp \overline{CD}$ with $|OA|=|OB|=|OC|=|OD|$.  
According Axiom (flat), $|AC|=|CB|=|BD|=|DA|=\sqrt{|OA|}$, so that $|\angle ACB|=|\angle CBD|=|\angle BDC|=|\angle DAC|=1/2$, or the geometric object $ACBD$ is a square.  

in Figure \ref{fig:squares} \textsl{Right}, when squares $ABCD$ and $CDEF$ are stacked along the side $\underline{CD}, C = \underline{*BF*} \perp \underline{*CD*}, D = \underline{*AE*} \perp \underline{*CD*}$, and the geometric object $ABFE$ is a rectangle. 
$|AB| = |CD|= |EF|$ is the distance between the two parallel straight lines $\underline{*AE*}$ and $\underline{*BF*}$.  
Thus, Theorem (flat) leads to the parallel postulate.

The stacking of identical squares can be continuous in both directions to form a measuring grid.
Both the side length and the orientation of the identical squares can be changed continuously, so that Euclidean space is homogeneous, isometric, isotropic, reflective, and similar, allowing both congruence relation and similarity relation between geometric objects.
In contrast, measuring grids can also be established in spherical space using equilateral triangles or polygons, but only for selected sizes \cite{Non-Euclidean Geometry}, so that a spherical space can have congruence relation but not similarity relation.

\subsection{Proof for Theorem (perpendicular)}

In addition to lemmas for Theorem (perpendicular), a few other involved important theorems for Euclidean surface geometry will also be presented.

\begin{thm}[triangle inequality]
$C \not \in \underline{AB} \iff |AB| < |AC| + |CB|$.
\end{thm}

A round geometric object separates the space into three parts:
\begin{defn}[inside or outside round] $\forall \underline{A:r}, B:$
\begin{itemize}
\item $B$ is \emph{inside} $\underline{A:r}$ if $|AB| < r$.
\item $B$ is \emph{outside} $\underline{A:r}$ if $|AB| > r$.
\end{itemize}
\end{defn}

If all points of a geometric object are inside/outside a round, the geometric object is inside/outside the round.

\begin{thm}[exclusively outside]
$\forall \underline{A:r_A}, \underline{B:r_B}: |AB| > r_A + r_B \iff \underline{A:r_A} $ and $ \underline{B:r_B}$ are outside of each other.
\end{thm}

\begin{thm}[inclusively outside]
$\forall \underline{A:r_A}, \underline{B:r_B}: |AB| = r_A + r_B \iff \underline{A:r_A} $ and $ \underline{B:r_B}$ are outside of each other except for
$\exists ! C: \{C\} = \underline{A:r_A} \cap \underline{B:r_B} \bigwedge C \in \underline{AB}$.
\end{thm}

\begin{thm}[exclusively inside]
$\forall \underline{A:r_A}, \underline{B:r_B}: |AB| < r_A - r_B \iff \underline{B:r_B} $ is inside $ \underline{A:r_A}$.
\end{thm}

\begin{thm}[inclusively outside]
$\forall \underline{A:r_A}, \underline{B:r_B}: |AB| = r_A - r_B \iff \underline{B:r_B} $ is inside $ \underline{A:r_A}$ except for
$\exists ! C: \{C\} = \underline{A:r_A} \cap \underline{B:r_B} \bigwedge B \in \underline{AC}$.
\end{thm}

As a further quantification of Axiom (strictly monotonic):
\begin{thm}[diminishing difference]
$D = C \perp \underline{*AB*} \bigwedge \forall E, F \in \underline{*AB*}:$
\begin{itemize}
\item $|DE| > |DF| \Longrightarrow 0 < |CE| - |DE| < |CF| - |DF| < |CD|$.

\item $|DE| \rightarrow 0 \Longrightarrow |CE| - |DE| \rightarrow |CD|$.

\item $|DE| \rightarrow +\infty \Longrightarrow |CE| - |DE| \rightarrow 0$.
\end{itemize}
\end{thm}
\begin{proof}
According to Figure \ref{fig:strictly_monotonic}:

$|CE| < |CF| + |EF| \Longrightarrow 0 < |CE| - |DE| < |CF| - |DF| < |CD|$.

$\exists ! G, H: \{G\} = \underline{E:|CE|} \cap \underline{*AB*} \;\bigwedge\; \{H\} = \underline{F:|CF|} \cap \underline{*AB*} \Longrightarrow G \in \underline{DH}$.

$|CE| \rightarrow +\infty \Longrightarrow |DG| = |CE| - |DE| \rightarrow 0$.
\end{proof}

\begin{figure}
\includegraphics{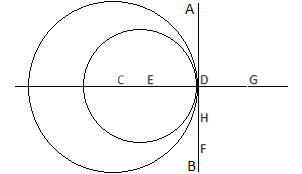}
\caption{Proof of Theorem (perpendicular).}
\label{fig:perpendicular}
\end{figure}

\begin{lem}[perpendicular 1]
$D = C \perp \underline{*AB*}: \forall E \in \underline{CD} \Longrightarrow D = E \perp \underline{*AB*}$
\end{lem}
\begin{proof}
According to Figure \ref{fig:perpendicular}, $\underline{E: |ED|}$  is inside $\underline{C: |CD|}$ except for $D$.
\end{proof}

\begin{lem}[perpendicular 2]
$D = C \perp \underline{*AB*}: \forall E \in \underline{*CD} \Longrightarrow D = E \perp \underline{*AB*}$
\end{lem}
\begin{proof}
Using Lemma (perpendicular 1) and proof by contradiction.
\end{proof}

\begin{lem}[perpendicular 3]
$D = C \perp \underline{*AB*}: \forall G: G \in \underline{*CD*} \bigwedge G \not \in \underline{*CD} \Longrightarrow D = G \perp \underline{*AB*}$
\end{lem}
\begin{proof}
According to Figure \ref{fig:perpendicular}, let $\forall E \in \underline{*CD}$ and $\exists ! H \in \underline{*AB*}: H = G \perp \underline{*CD*} \Longrightarrow |ED| < |EH| \bigwedge |GH| < |DG|$.

$|EG| < |EH| + |GH| \Longrightarrow 0 < |DG| - |GH| < |EH| - |DE|$.  

$|DE| \rightarrow +\infty \Longrightarrow |EH| - |DE| \rightarrow 0 \Longrightarrow |DG| \rightarrow |DH| \Longrightarrow H = D$.
\end{proof}

The combination of Lemma (perpendicular 2) and (perpendicular 3) suggests $D = \underline{*CD*} \perp \underline{*AB*}$.

It is desired that Lemma (perpendicular 3) be proven by Theorem (inclusively outside), rather than by Theorem (diminishing difference) which is no longer be local.

\begin{lem}[perpendicular 4]
$D = \underline{*CD*} \perp \underline{*AB*} \Longrightarrow D = \underline{*AB*} \perp \underline{*CD*}$.
\end{lem}
\begin{proof}
According to Figure \ref{fig:perpendicular}, let $H = \underline{*AB*} \perp \underline{*CD*}$.
  
$\forall F \in \underline{*AB*} \Longrightarrow \exists ! E \ne G: E, G \in \underline{*CD*}: |DF| = |DE| = |DG| \Longrightarrow |DF| < |EF|, |FG| \longrightarrow H \in \underline{EG}$.

$|DF| \rightarrow 0 \Longrightarrow |EG| \rightarrow 0 \Longrightarrow H = D$.
\end{proof}

Another way to prove Lemma (perpendicular 4) is from Axiom (measurable) which indicates that two straight lines intersect at most one point.

\subsection{Proof for Theorem (flat)}

Most of this section states properties of an isometric and isotropic surface.

\begin{defn}[triangle]
When $C \not \in \underline{*AB*}$, $\underline{AB} \cup \underline{BC} \cup \underline{CA}$ is defined as a \emph{triangle} $\triangle ABC$, in which each of $\underline{AB}$, $\underline{BC}$, $\underline{CA}$ is defined as a \emph{side}, and each of $\angle ABC$, $\angle BCA$, and $\angle CAB$ is an inner angle.
\end{defn}

\begin{defn}[equivalent triangles]
If two triangles $\forall \triangle ABC, \triangle DEF$ have pair-wise equal length sides, and equal inner angles by corresponding pairs, these two triangles are \emph{equivalent triangles} $\triangle ABC = \triangle DEF$.
\end{defn}

\begin{defn}[side-angle-side equivalence]
The relation $\forall \triangle ABC, \triangle DEF: |AB| = |DE| \bigwedge |AC| = |DF| \bigwedge |\angle CAB|_{AB} = |\angle FDE|_{|DE|} \Longrightarrow \triangle ABC = \triangle DEF$ is defined as \emph{side-angle-side equivalence}.
\end{defn}
The side-angle-side equivalence \cite{Plane} leads to other equivalence such as side-side-side equivalence or side-angle-angle equivalence \cite{Plane}, which will be used directly in the following.

\begin{defn}[isometric and isotropic]
A surface is \emph{isometric and isotropic} if side-angle-side equivalence holds.
\end{defn}

A measurement of the angle distance ratio is an application of the side-angle-side equivalence:
\begin{lem}[flat 1]
On a surface, if the right ratio is a constant, the surface is isometric and isotropic.
\end{lem}

\begin{figure}
\includegraphics{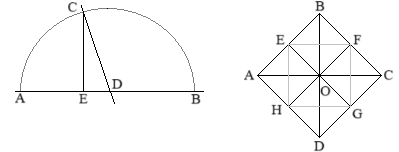}
\caption{Proof of Theorem (flat).}
\label{fig:flat}
\end{figure}

\begin{lem}[flat 2]
$\forall r > 0 \bigwedge \forall C \not \in \underline{*AB*} \bigwedge \forall D \in \underline{*AB*}: |\angle ADC|_r = |\angle BDC|_r 
\Longrightarrow D = C \perp \underline{*AB*}$.
\end{lem}
\begin{proof}
In Figure \ref{fig:flat} \emph{left}, $A, B, C \in \underline{D:r} \bigwedge |AC| = |BC|$.
Let $E = C \perp \underline{*AB*}$ and assume $E \in \underline{AD}$.
According to Axiom (strictly monotonic), $|AC| < |BC|$ unless $D = E$.
\end{proof}

In Figure \ref{fig:flat} \emph{right}, construct $O = \overline{AC} \perp \overline{BD}$ with $|OA|=|OB|=|OC|=|OD|=1$.
Let $E, F, G$, and $H$ be the midpoints of $\overline{AB}, \overline{BC}, \overline{CD}$, and $\overline{DA}$, respectively.
Let $|AB|=|BC|=|CD|=|DA|=\alpha, |OE| = \beta$, so that the right ratio is $\alpha^2 / 4$.
\begin{enumerate}
\item $|\angle DAB|_{\alpha}=|\angle ABC|_{\alpha}=|\angle BCD|_{\alpha}=|\angle CDA|_{\alpha} = 1/ \alpha^2$.

\item $\triangle EBF = \triangle FCG \Longrightarrow |EF| = |FG| = |GH| = |HE|$.

\item $\triangle OEA = \triangle OEB 
\Longrightarrow E = \underline{OE} \perp \underline{AB} \bigwedge F = \underline{OF} \perp \underline{BC} \bigwedge G = \underline{OG} \perp \underline{CD} \bigwedge H = \underline{OH} \perp \underline{DA}$.

\item $\triangle OEA = \triangle OFC \Longrightarrow |OE| = |OF| = |OG| = |OH| = \beta$.

\item $\triangle OEF = \triangle OFG \Longrightarrow O = \underline{EG} \perp \underline{FH} \Longrightarrow |EF| = |FG| = |GH| = |HE| = \alpha \beta$. 

\item $\|\angle HEF|_{\frac{\alpha \beta}{2}} = |\angle EFG|_{\frac{\alpha \beta}{2}} = |\angle FGH|_{\frac{\alpha \beta}{2}}= |\angle GHE|_{\frac{\alpha \beta}{2}} = 1 /\alpha^2 
\Longrightarrow$  an angle distance ratio does not depends on the angle measuring distances, so that it is a constant, which is consistent with the constant right ratios.

\item $|\angle HEF| = |\angle DAB| 
\Longrightarrow |EF|= 1 
\Longrightarrow |OE| = \beta = 1 /\alpha$.

\item $|\angle AEO| = \alpha^2 / 4 = |AO|^2 / |OE|^2 /4 \Longrightarrow |OE| = |AE| \Longrightarrow 1 /\alpha = \alpha /2 \Longrightarrow $ the right ratio is $1/2$.
\end{enumerate}

\section{Beyond Euclidean Geometry}

\subsection{Right Ratio}

\begin{figure}
\includegraphics[scale=0.6]{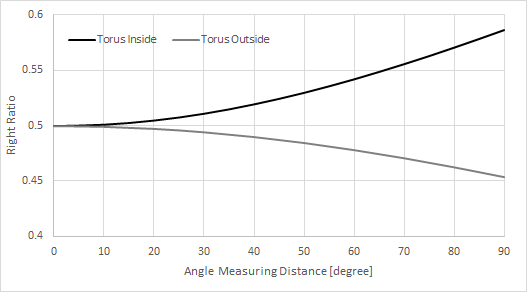}
\caption{
The right ratios at the centers of the inside surface and the outside surface of a torus \cite{Non-Euclidean Geometry}.  
The x-axis is the angle measuring distance normalized as the angle from the center crossing the torus.  
The torus has an inner radius $2$ and an outer radius $3$. 
}
\label{fig:right ratios}
\end{figure}

Right ratios correlate with corresponding Gaussian curvatures. 
For example, a torus \cite{Non-Euclidean Geometry} has negative Gaussian curvatures at its inside surface, and positive Gaussian curvatures at its outside surface. 
Figure \ref{fig:right ratios} shows that the right ratios at the center of the inside surface of a torus are more than 1/2 and increase with the angle measuring distance, while the right ratios at the center of the outside surface of the torus are less than 1/2 and decrease with the angle measuring distance.
Figure \ref{fig:right ratios} also shows that while the first derivative of the local right ratio is zero, the second derivative correlates to the value of Gaussian curvature. 
While Gaussian curvature is an explicit measurement of surface curvature, the local right ratio is an implicit measurement.  
The local right ratio simpler than Riemann curvature tensor \cite{Differential Geometry} for local measurement.

As a conclusion of Theorem (flat), when the first derivative of the local right ratio is zero, the surface is manifold \cite{Differential Geometry} at the point.  
Figure \ref{fig:right ratios} presents such an example.

On a conical surface, the right ratios are not continuous when the tip is included or excluded in the corresponding right angles, so that the continuity of the right ratios excludes such singular points.  
Such discontinuity holds for local right ratio as well.
\begin{cnj}[manifold]
If the local right ratios are smooth on a surface, the surface is a manifold at each point.
\end{cnj}

\subsection{Isometric and Isotopic}

It is desirable that an axiom set distinguishes Euclidean surfaces from other surfaces only in the last axiom.
However, such desire conflicts with the original Euclidean axiom set: on a spherical surface, between two opposing poles $A$ and $B$ there are infinitive $\underline{AB}$, and between a point and a geodesic which is not the equatorial of the point there are two tangent points.
An alternative is to use isometric and isotopic as the foundation for both Euclidean and spherical surfaces \cite{Isometric}, which has not been proven rigorously.
\begin{cnj}[isometric and isotopic]
If a 2-sided surface is isometric and isotopic, it is either Euclidean or spherical.
\end{cnj}
On an isometric and isotopic surface, the angle distance ratios depends on the angle measuring distances only, which may be used to prove Conjecture (isometric and isotopic) rigorously.

\section{Discussion}

\subsection{Summary}

With Axiom (strictly monotonic) added, and the space being limited to topologically simplest 2-dimensional space, the modernized Euclidean axiom set can describe Euclidean surface accurately and completely, to be:
\begin{enumerate}
\item  Axiom (measurable) is equivalent to Euclidean Axiom 3.
\item  Axiom (continuous) is equivalent to Euclidean Axiom 1.
\item  Axiom (boundless) is equivalent to Euclidean Axiom 2.
\item  Axiom (strictly monotonic) is a new axiom.
\item  Definition (surface) limits Euclidean to be topologically simplest 2-dimensional space.
\item Theorem (flat) is equivalent to Euclidean Axiom 4.
\end{enumerate}
Axiom (measurable), Axiom (continuous), and Axiom (boundless) can be combined as one axiom: $\forall A, B \Longrightarrow \exists ! \underline{*AB*}$.  
But such simplification hides the progressive dependency between the three axioms.

As the new and simpler surface curvature measurement analytic, right ratio can replace both Gaussian curvature and Riemann curvature tensor.

\subsection{Improving the New Axiom Set}

From intuition, Theorem (flat) may be improved as:
\begin{cnj}[flat]
If on a surface, all right ratios are a constant at one point, they are $1/2$.
\end{cnj}

Axiom (strictly monotonic) is a new axiom that seems to be implied in all the existing axiom sets for Euclidean geometry.
It needs further scrutiny.
When a point is very close to a geodesic, it always holds.
If the condition for Axiom (strictly monotonic) can be quantified, it may also be reduced to a theorem as how the 4th Euclidean axiom is reduced to Theorem (flat).

\subsection{Application to Measure Space Curvature for Astrophysics}

Measuring space curvature experimentally is both important and active research today in astrophysics, such as the detection of gravitational waves \cite{Gravitational Waves}. 
All current methods are based on the parallel postulate, such as the sum of angles of a triangle, so they can only measure integral effects of space curvature over long distance \cite{Space Curvature}. 
Right ratios provide space curvature measurement at or near each point, so it could be an important improvement.  
For example, in Figure \ref{fig:squares} Left, if the length difference between each pair of arms of three Michelson interferometers (as what is used in the detection of gravitational waves) $OA-OB$, $OB-OC$, and $BA-BC$ is all within one wavelength, $O = \underline{CD} \perp \underline{AB}$ with $|OA| = |OB| = |OC|$, and $C = \underline{AC} \perp \underline{BC}$ with $|AC| = |BC|$ are achieved accurately.  
Figure \ref{fig:squares} left contains measurements for 8 right ratios, to accurately measure the space curvature according to Figure \ref{fig:right ratios}.
Such a setup is an improvement over the current method of using one Michelson interferometer at one site to measure the violation of Axiom (measurable) only.

\ifdefined\AUTHOR

\subsection{Educational Purpose}

In US, the public education on mathematics is not satisfactory to prepare students for science and engineering. 
Best students learn mathematics by themselves, such as through the requirements of AMC, AIME and USAMO competitions, which are all way ahead of the normal curriculum. 
For example, Euclidean geometry appears in AMC8 for 8th graders as choice selections, but all courses related to geometries including 2D and 3D Euclidean geometry, analytic geometry, and trigonometry are crammed in the 11th grade for Alice (who is the 2nd author) officially. 
The self-learning for completions or the tight curriculum schedule only allows a quick run of basic knowledge of each subject matter, and only a touch of mathematical proving. 
Both approaches miss the chance for a new student to ask naive but sometimes profound questions, which is the major source of creativeness and understanding in mathematics. 
The tight schedule also reduces teachers from inspirers to instructors. AP courses or other advanced courses make the trend of not understanding a subject matter thoroughly worse. 
In contrast, this paper is generated in the learning process of both authors over 2 years when Alice learned surface Euclidean geometry for the first time in a relaxed manner from Chengpu, when the naive but sometimes deep questions from Alice broke the for-granted understanding by Chengpu, starting from the question of how to measure an angle.  
Relaxed time also allows both authors to use mathematical proofs as a reasoning process. 
Both authors hope that the mathematical teaching in US public schools should be improved on the scheduling and the emphasis of mathematical curriculum.

\section{Acknowledgment}

Chengpu Wang feels deeply grateful for:
\begin{itemize}
\item the wonderful teaching by Mr. Jianye Liu from The Middle School of Peking University, and

\item the unique teaching for challenging status quo of knowledge, by Dr Paul Hough from Brookhaven National Lab, including the courage to reinvent wheels, either as a learning process or as an inventory method.
\end{itemize}

The first attempt resulted in a fatally flawed self-published paper How to Define a Flat Plane (Kindle, 2017). Mr. Victor Aguilar, the author of Geometry-Do (2019), spent time reading the first paper and found the flaw. In 2023, Dr. Oliver Attie, a mathematician on differential geometry and bioinformatics, reviewed the paper carefully, and provided many valuable suggestions. Many personal friends help with proofreading the draft of the paper.

\fi

% ------------------------------------------------------------------------
\end{document}